\documentclass[12pt]{amsart}
\usepackage{amssymb,bbold}
\usepackage[all]{xy}

\theoremstyle{plain}
\newtheorem{theorem}{Theorem}

\newtheorem{lemma}[theorem]{Lemma}
\newtheorem{proposition}[theorem]{Proposition}

\theoremstyle{definition}

\newtheorem{remark}[theorem]{Remark}
\newtheorem{example}[theorem]{Example}
\newtheorem{question}{Question}

\def\one{\mathbb 1}

\begin{document}
\baselineskip 18pt

\title[Fremlin tensor product of Banach lattices ]
      {Fremlin tensor product of Banach lattices and unbounded convergences}
\author[O.~Zabeti]{Omid Zabeti}
\address[O.~Zabeti]
  {Department of Mathematics, Faculty of Mathematics, Statistics, and Computer science,
   University of Sistan and Baluchestan, Zahedan,
   P.O. Box 98135-674. Iran}
\email{o.zabeti@gmail.com}
\maketitle

\begin{abstract}
Let $E$ and $F$ be Banach lattices. In this short note, we consider several conditions under them, the Fremlin projective tensor product $E\widehat{\otimes}F$ respects either the unbounded norm convergence or unbounded absolute weak convergence. 

\end{abstract}

\date{\today}

\maketitle
\section{Motivation and Preliminaries}
Assume that $E$ and $F$ are vector lattices. In \cite{Fremlin:72}, Fremlin constructed the vector lattice
$E\overline{\otimes}F$, named as the Fremlin tensor product of $E$ and $F$. In this structure, the
algebraic tensor product $E\otimes F$ can be considered as an ordered vector subspace of
$E\overline{\otimes}F$.

Now, let $E$ and $F$ be Banach lattices, Fremlin  in \cite{Fremlin:74} provided the Fremlin
projective tensor product $E\widehat{\otimes}F$, which is itself a Banach lattice. In this setting, the algebraic tensor product $E\otimes F$ is norm dense in
$E\widehat{\otimes}F$. Furthermore, 
$E\overline{\otimes}F$ can be considered as a norm-dense vector sublattice of
$E\widehat{\otimes}F$; see \cite[1A]{Fremlin:74} for a comprehensive details.

In \cite{Z:25}, it was shown that the Fremlin tensor product of Archimedean vector lattices preserves unbounded order convergence as well as order convergence. On the other hand, in a  Banach lattice, two natural notions of unbounded convergence play a central
role: unbounded norm convergence and unbounded absolute weak convergence. This leads to the
following question.

\begin{question}\label{0000}
Does the Fremlin projective tensor product of two Banach lattices preserve either unbounded norm
convergence or unbounded absolute weak convergence?
\end{question}

 In the present
paper, we  establish, under mild assumptions, that the Fremlin projective tensor product between Banach lattices
preserves both unbounded norm convergence and unbounded absolute weak convergence. 

We begin by recalling several notions of unbounded convergence in Banach lattices. 
Let $E$ be a Banach lattice and let $(x_\alpha)\subseteq E$ be a net. 
The net $(x_\alpha)$ is said to be \emph{unbounded order convergent} (uo-convergent) to $x\in E$ if
\[
|x_\alpha-x|\wedge u \xrightarrow{o} 0
\quad\text{for every } u\in E_+;
\]
in the other words, the net $(|x_{\alpha}-x|\wedge u)$ is order convergent to zero. 

It is called \emph{unbounded norm convergent} (un-convergent) to $x\in E$ if
\[
\bigl\|\,|x_\alpha-x|\wedge u\,\bigr\|\to 0
\quad\text{for all } u\in E_+.
\]
Furthermore, $(x_\alpha)$ is said to be \emph{unbounded absolutely weakly convergent} to $x\in E$, denoted by
$x_\alpha\xrightarrow{uaw}x$, if
\[
|x_\alpha-x|\wedge u \xrightarrow{w} 0
\quad\text{for every } u\in E_+;
\]
that is the net $(|x_{\alpha}-x|\wedge u)$ is weakly null.

For further information on these notions and related concepts, we refer to
\cite{Den:17, GTX:17, KMT, Z:18}. 
Moreover, for general terminology and background on vector and Banach lattices, see
\cite{AB1, AB}.

Next, we briefly recall some facts concerning the Fremlin tensor product of vector
and Banach lattices; see \cite{Fremlin:72, Fremlin:74} for details.

Suppose that $E$ and $F$ are Archimedean vector lattices. 
Fremlin in \cite{Fremlin:72} established  a tensor product $E\overline{\otimes}F$,
which is again an Archimedean vector lattice containing the algebraic tensor product
$E\otimes F$ as a vector subspace.

Now let $E$ and $F$ be Banach lattices, Fremlin in \cite{Fremlin:74} further constructed the
projective tensor product $E\widehat{\otimes}F$, which is a Banach lattice. 
This space is the completion of $E\otimes F$ with respect to the projective norm
$\|\cdot\|_{|\pi|}$, which is a lattice cross norm satisfying
\[
\|x\otimes y\|_{|\pi|}=\|x\|\,\|y\|
\quad\text{for all } x\in E,\; y\in F.
\]
In this structure, $E\overline{\otimes}F$ is a norm-dense vector sublattice
of $E\widehat{\otimes}F$.
\section{main results}

In \cite{Z:25}, it is proved that the Fremlin tensor product of vector lattices preserves unbounded order convergence. More precisely, we state the main result of that work (see \cite[Theorem 7]{Z:25}).
\begin{theorem}
Suppose $E$ and $F$ are Archimedean vector lattices. Moreover, assume that $x_{\alpha}\xrightarrow{uo}x$ in $E$ and $y_{\beta}\xrightarrow{uo}y$ in $F$. Then, $x_{\alpha}\otimes y_{\beta}\xrightarrow{uo} x\otimes y$ in the Fremlin tensor product $E\overline{\otimes}F$.
\end{theorem}
On the other hand, as mentioned before, we have the Fremlin projectve tensor product between Bananch lattices.
 In what follows, we investigate some mild assumptions, under which, Fremlin projective tensor product between Banach lattices preserves unbounded norm convergence as well as unbounded absolute weak convergence. Note that both convergences are topological. In other words, any Banach lattice with the assumed convergence, forms a locally solid vector lattice. 

\begin{lemma}\label{304}
Suppose $E$ and $F$ are Banach lattices such that the Fremlin projective tensor product $E\widehat{\otimes} F$ is order continuous. Then, both $E$ and $F$ are order continuous.  
\end{lemma}
\begin{proof}
We show that $E$ is order continuous; the proof for $F$ is analogous.
Suppose, to the contrary, that $E$ is not order continuous.
Then there exists a disjoint, order-bounded sequence $(x_n)\subseteq E$
that is not norm null.
Choose any nonzero positive element $y\in F$.
It is easy to see that the sequence $(x_n\otimes y)$ is disjoint and
order bounded in $E\widehat{\otimes} F$.
By the assumption, this sequence is norm null.
However, since the projective tensor norm is a cross norm,
it follows that $(x_n)$ must be norm null, which is a contradiction.

\end{proof}
Observe that the converse of Lemma~\ref{304} is not true,
even for reflexive Banach lattices,
as shown in \cite[Example~4C]{Fremlin:74}.
More precisely, let $E=F=L^2[0,1]$.
Then the Fremlin projective tensor product $E\widehat{\otimes}F$
is not order complete, and hence the projective norm on
$E\widehat{\otimes}F$ is not order continuous,
although both $E$ and $F$ are order continuous (indeed, reflexive)
Banach lattices.
\begin{proposition}\label{202}
Suppose $E$ and $F$ are Banach lattices. Assume that $(x_{\alpha})\subseteq E$ is $un$-null and $(y_{\beta})\subseteq F$ is eventually norm bounded. Then, the double net $(x_{\alpha}\otimes y_{\beta})$ is $un$-null in $E\widehat{\otimes}F$.
\end{proposition}
\begin{proof}
By the assumption, $\|y_{\beta}\|\leq 1$ for sufficiently large $\beta$. Suppose $w\in (E\overline{\otimes}F)_{+}$. By  \cite[1A (d)]{Fremlin:74}), there exist $x_0\in E_{+}$ and $y_0\in F_{+}$ with $w\leq x_0\otimes y_0$. Now, we use the following inequality.
\[|x_{\alpha}\otimes y_{\beta}|\wedge w\leq (|x_{\alpha}|\otimes |y_{\beta}|)\wedge (x_0\otimes y_0)\]
\[\leq (|x_{\alpha}|\wedge x_0)\otimes (|y_{\beta}|\vee y_0).\]
Therefore,
\[\||x_{\alpha}\otimes y_{\beta}|\wedge w\|\leq \||x_{\alpha}|\wedge x_0\|\||y_{\beta}|\vee y_0\|\leq \||x_{\alpha}|\wedge x_0\|(\|y_0\|+1)\rightarrow 0.\]
We conclude that $x_{\alpha}\otimes y_{\beta}\xrightarrow{un}0$ in  $E\overline{\otimes} F$ and so that in $E\widehat{\otimes}F$ by \cite[Theorem 4.3]{KMT}; see also \cite[Lemma 3.5]{Tay:19}.
\end{proof}
\begin{theorem}\label{2222}
Suppose $E$ and $F$ are Banach lattices such that $E\widehat{\otimes}F$ is order continuous. If $(x_{\alpha})\subseteq E$ is $un$-convergent to $x\in E$ and $(y_{\beta})\subseteq F$ is $un$-convergent to $y\in F$, then the double net $(x_{\alpha}\otimes y_{\beta})$ is $un$-convergent to $x\otimes y$ in $E\widehat{\otimes}F$. 
\end{theorem}
\begin{proof}
First, note that by Lemma~\ref{304}, both $E$ and $F$ are order continuous.

By \cite[Theorem~3.2]{Den:17}, there exists a subsequence
$(\alpha_k)$ of indices and a disjoint sequence $(d_k)$ in $E$
such that
\[
\|x_{\alpha_k}-x-d_k\|\longrightarrow 0.
\]
Moreover, by \cite[Corollary 3.5]{Den:17}, we can find an increasing sequence $(\beta_k)$ of indices such that 
\[y_{\beta_k}-y\xrightarrow{un}0.
\]
By Proposition~\ref{202},
\[
(x_{\alpha_k}-x-d_k)\otimes (y_{\beta_k}-y)\xrightarrow{un}0,
\]
\[
(x_{\alpha_k}-x)\otimes y\xrightarrow{un}0,
\qquad
x\otimes (y_{\beta_k}-y)\xrightarrow{un}0.
\]
Furthermore, observe that the sequence
$d_k\otimes (y_{\beta_k}-y)$ is disjoint in $E\widehat{\otimes}F$,
and hence $un$-null by \cite[Proposition~3.5]{KMT}.
Therefore, we have
\begin{align*}
x_{\alpha_k}\otimes y_{\beta_k}-x\otimes y
&=(x_{\alpha_k}-x-d_k)\otimes (y_{\beta_k}-y)
+ d_k\otimes (y_{\beta_k}-y) \\
&\quad + x\otimes (y_{\beta_k}-y)
+ (x_{\alpha_k}-x)\otimes y
\xrightarrow{un}0.
\end{align*}

This argument applies to every subnet of $(x_{\alpha})$ and
$(y_{\beta})$.
Since $un$-convergence is topological, we conclude that
\[
x_{\alpha}\otimes y_{\beta}\xrightarrow{un} x\otimes y.
\]

\end{proof}
\begin{question}\label{90}
Is the order continuity of the Fremlin projective tensor product $E\widehat{\otimes}F$ a necessary hypothesis in Theorem \ref{2222}?
\end{question}
We give a negative answer to this question by the following example.
\begin{example}
Assume that $K_1$ and $K_2$ are compact Hausdorff spaces.
It is well known that
\[
C(K_1)\widehat{\otimes} C(K_2)=C(K_1\times K_2)
\]
(see \cite{Fremlin:74} for more details).
Moreover, by \cite[Theorem~2.3]{KMT}, $un$-convergence in $C(K)$-spaces
coincides with norm convergence.

Now assume that $(f_{\alpha})\subseteq C(K_1)$ is $un$-null and that
$(g_{\beta})\subseteq C(K_2)$ is also $un$-null.
Then both nets are norm null.
On the other hand,
\[
(f_{\alpha}\otimes g_{\beta})(t,s)=f_{\alpha}(t)g_{\beta}(s),
\quad t\in K_1,\ s\in K_2.
\]
Moreover,
\[
\|f_{\alpha}\otimes g_{\beta}\|
=\|f_{\alpha}\|\,\|g_{\beta}\|\longrightarrow 0.
\]
Therefore, the double net $(f_{\alpha}\otimes g_{\beta})$
is $un$-null in $C(K_1\times K_2)$.

\end{example}
We proceed to study the relationship between unbounded absolute weak convergence and the Fremlin projective tensor product.
\begin{proposition}\label{2022}
Let $E$ and $F$ be Banach lattices. Assume that $(x_{\alpha})\subseteq E$ is
$uaw$-null and that $(y_{\beta})\subseteq F$ is eventually order bounded.
Then the double net $(x_{\alpha}\otimes y_{\beta})$ is $uaw$-null in
$E\widehat{\otimes}F$.
\end{proposition}

\begin{proof}
By assumption, there exists $y_{1}\in F_{+}$ such that
\[
|y_{\beta}|\leq y_{1}
\quad \text{for all sufficiently large } \beta.
\]
Let $w\in (E\overline{\otimes}F)_{+}$.
By \cite[1A(d)]{Fremlin:74}, there exist $x_{0}\in E_{+}$ and
$y_{0}\in F_{+}$ such that
\[
w\leq x_{0}\otimes y_{0}.
\]

Fix
\[
\widetilde{f}\in (E\widehat{\otimes}F)'_{+}
= B^{r}(E\times F)_{+},
\]
where $B^{r}(E\times F)$ denotes the Banach lattice of all regular bounded
bilinear forms on $E\times F$; see \cite[1I]{Fremlin:74}.
Under this identification, there exists a unique
$f\in B^{r}(E\times F)_{+}$ such that
\[
\widetilde{f}(x\otimes y)=f(x,y),
\qquad x\in E,\ y\in F.
\]

We now estimate
\begin{align*}
\widetilde{f}(|x_{\alpha}\otimes y_{\beta}|\wedge w)
&\leq \widetilde{f}((|x_{\alpha}|\otimes |y_{\beta}|)\wedge (x_{0}\otimes y_{0}))\\
&\leq \widetilde{f}\big((|x_{\alpha}|\wedge x_{0})\otimes (|y_{\beta}|\vee y_{0})\big)\\
&\leq \widetilde{f}\big((|x_{\alpha}|\wedge x_{0})\otimes (y_{1}\vee y_{0})\big)\\
&= f(|x_{\alpha}|\wedge x_{0},\, y_{1}\vee y_{0}).
\end{align*}

Define the mapping
\[
f_{\,y_{1}\vee y_{0}}\colon E\to \mathbb{R},
\qquad
f_{\,y_{1}\vee y_{0}}(x)=f(x,\,y_{1}\vee y_{0}).
\]
Then $f_{\,y_{1}\vee y_{0}}$ is a positive linear functional on $E$.
Since $(x_{\alpha})$ is $uaw$-null in $E$, we obtain
\[
f(|x_{\alpha}|\wedge x_{0},\, y_{1}\vee y_{0})\longrightarrow 0.
\]
Hence,
\[
\widetilde{f}(|x_{\alpha}\otimes y_{\beta}|\wedge w)\longrightarrow 0.
\]

Therefore,
$x_{\alpha}\otimes y_{\beta}\xrightarrow{uaw}0$
in $E\overline{\otimes}F$, and consequently also in
$E\widehat{\otimes}F$ by \cite[Lemma~3.4]{Tay:19}.
\end{proof}

\begin{theorem}\label{2023}
Let $E$ and $F$ be order continuous Banach lattices.
If $(x_{\alpha})\subseteq E$ is $uaw$-convergent to $x\in E$ and
$(y_{\beta})\subseteq F$ is $uaw$-convergent to $y\in F$, then the double net
$(x_{\alpha}\otimes y_{\beta})$ is $uaw$-convergent to $x\otimes y$ in
$E\widehat{\otimes}F$.
\end{theorem}

\begin{proof}
By \cite[Theorem~4]{Z:18}, we have $x_{\alpha}\xrightarrow{un}x$ and $y_{\beta}\xrightarrow{un}y$.
Hence, by \cite[Theorem~3.2]{Den:17}, there exist an increasing sequence
$(\alpha_k)$ of indices and a disjoint sequence $(d_k)$ in $E$ such that
\[
\|x_{\alpha_k}-x-d_k\|\longrightarrow 0.
\]
Passing to a further subsequence if necessary, we may assume that
$(x_{\alpha_k}-x-d_k)$ is order bounded.

Furthermore, by \cite[Corollary 3.5]{Den:17} and also \cite[Theorem 4]{Z:18}, there exists an increasing sequence $(\beta_k)$ of indices such that 
\[
y_{\beta_k}-y\xrightarrow{un}0.
\]
By Proposition~\ref{2022}, it follows that
\[
(x_{\alpha_k}-x-d_k)\otimes (y_{\beta_k}-y)\xrightarrow{uaw}0
\]
and also that
\[
(x_{\alpha_k}-x)\otimes y\xrightarrow{uaw}0
\quad\text{and}\quad
x\otimes (y_{\beta_k}-y)\xrightarrow{uaw}0.
\]

Moreover, since $(d_k)$ is a disjoint sequence in $E$, the sequence
$(d_k\otimes (y_{\beta_k}-y))$ is disjoint in $E\widehat{\otimes}F$ and hence
$uaw$-null by \cite[Lemma~2]{Z:18}.

Combining these observations, we obtain
\begin{align*}
x_{\alpha_k}\otimes y_{\beta_k}-x\otimes y
&=(x_{\alpha_k}-x-d_k)\otimes (y_{\beta_k}-y)
+d_k\otimes (y_{\beta_k}-y)\\
&\quad +x\otimes (y_{\beta_k}-y)
+(x_{\alpha_k}-x)\otimes y
\xrightarrow{uaw}0.
\end{align*}

Since this argument applies to every subnet of $(x_{\alpha})$ and $(y_{\beta})$, we conclude that
\[
x_{\alpha}\otimes y_{\beta}\xrightarrow{uaw}x\otimes y
\quad\text{in } E\widehat{\otimes}F.
\]
\end{proof}

We are thus led to the following question.

\begin{question}\label{900}
Is the order continuity of the components a necessary assumption
in Theorem~\ref{2023}?
\end{question}

We now provide a negative answer to this question.

\begin{lemma}\label{1}
Let $K$ be a compact Hausdorff space and let $(f_n)\subseteq C(K)$.
Then $f_n\xrightarrow{uaw}0$ if and only if $f_n\to 0$ pointwise on $K$.
\end{lemma}

\begin{proof}
Assume first that $f_n\xrightarrow{uaw}0$.
Then $|f_n|\wedge \one \xrightarrow{w}0$.
By the Riesz representation theorem (see \cite[p.~203]{AB} for details),
it follows that
\[
|f_n(t)|\longrightarrow 0
\quad \text{for each } t\in K,
\]
and hence $f_n(t)\to 0$ pointwise on $K$.

Conversely, assume that $f_n\to 0$ pointwise on $K$.
Then $|f_n|\wedge \one \to 0$ pointwise.
Since the sequence $(|f_n|\wedge \one)$ is norm bounded,
the Riesz representation theorem again yields
\[
|f_n|\wedge \one \xrightarrow{w}0.
\]
Therefore, by \cite[Lemma~3]{Z:18}, we conclude that
$f_n\xrightarrow{uaw}0$.
\end{proof}

We now obtain the following result.

\begin{theorem}
Let $K_{1}$ and $K_{2}$ be compact Hausdorff spaces, and let
$(f_n)\subseteq C(K_{1})$ and $(g_n)\subseteq C(K_{2})$ be $uaw$-null
sequences. Then the sequence $(f_n\otimes g_n)$ is $uaw$-null in
$C(K_{1}\times K_{2})$.
\end{theorem}

\begin{proof}
By \cite[3E]{Fremlin:74}, we have the lattice isometric identification
\[
C(K_{1})\widehat{\otimes} C(K_{2}) = C(K_{1}\times K_{2}).
\]
Moreover,
\[
(f_n\otimes g_n)(t,s)=f_n(t)g_n(s)
\qquad
(t,s)\in K_{1}\times K_{2}.
\]

By Lemma~\ref{1}, the sequences $(f_n)$ and $(g_n)$ converge pointwise to
zero on $K_{1}$ and $K_{2}$, respectively.
Consequently,
\[
f_n\otimes g_n \longrightarrow 0
\quad\text{pointwise on } K_{1}\times K_{2}.
\]
Applying Lemma~\ref{1} once more, we conclude that
\[
f_n\otimes g_n \xrightarrow{uaw}0
\quad\text{in } C(K_{1}\times K_{2}).
\]
\end{proof}

In Theorem~\ref{2023}, we showed that for Banach lattices $E$ and $F$, the order
continuity of  the both factors is sufficient to guarantee that the Fremlin
projective tensor product preserves unbounded absolute weak convergence.
We now seek alternative conditions that yield the same conclusion.


\begin{remark}
Observe that, in general, there is no direct relationship between order
continuity of a Banach lattice $E$ and separability of its dual $E'$.
For example, $\ell_{1}$ is order continuous, whereas its dual
$\ell_{\infty}$ is not separable.
Conversely, the Banach lattice $c$ of convergent sequences is not order
continuous, while its dual $\ell_{1}$ is separable.
\end{remark}

It was shown in \cite[Theorem~3.2]{KMT} that the $un$-topology on a Banach lattice
$E$ is metrizable if and only if $E$ admits a quasi-interior point.
Moreover, by \cite[Theorem~4]{Z:18}, if $E$ is order continuous, the same
characterization holds for the $uaw$-topology.
However, order continuity is essential in this result and cannot be omitted,
even when the Banach lattice has a strong unit.

Indeed, consider $E=\ell_{\infty}$.
It is easy to see that the $uaw$-topology coincides with the weak topology on the
closed unit ball of $E$ (see \cite[Theorem~7]{Z:18}).
Consequently, the $uaw$-topology on $E$ is not metrizable, since $E'$ is not
separable.
Nevertheless, separability remains a useful assumption in this context.

\begin{lemma}\label{0003}
Suppose $E$ is a Banach lattice whose dual space $E'$ is separable.
Then the $uaw$-topology on $E$ is metrizable.
\end{lemma}

\begin{proof}
Since $E'$ is separable, it follows that $E$ itself is separable.
In particular, $E$ admits a quasi-interior point $e$.
Let $(f_k)_{k\in\mathbb N}$ be a countable dense subset of $E'$.

Define a mapping $d \colon E\times E \to \mathbb R$ by
\[
d(x,y)
=
\sum_{k=1}^{\infty}
\frac{|f_k(|x-y|\wedge e)|}{1+|f_k(|x-y|\wedge e)|}\,
\frac{1}{2^k}.
\]
It is straightforward to verify that $d$ is a translation-invariant metric on
$E$.

Moreover, for a net $(x_\alpha)$ in $E$, we have $d(x_\alpha,x)\to 0$ if and only
if
\[
f_k(|x_\alpha-x|\wedge e)\to 0
\quad \text{for all } k\in\mathbb N,
\]
which is equivalent to
\[
f(|x_\alpha-x|\wedge e)\to 0
\quad \text{for every } f\in E',
\]
by density of $(f_k)$ in $E'$.
The latter condition is precisely $uaw$-convergence.
Hence, the metric $d$ generates the $uaw$-topology on $E$.
For further details, see \cite[Lemma~3]{Z:18}.
\end{proof}

\begin{remark}
Note that if a locally solid vector lattice is metrizable by a
translation-invariant metric, then its topological completion is also a
metrizable locally solid vector lattice.
Moreover, the completion admits a translation-invariant metric that
extends the original one.
\end{remark}

\begin{theorem}
Assume that $E$ and $F$ are Banach lattices such that $E'$ is separable and $F$ is order continuous.
If $(x_{\alpha})\subseteq E$ is $uaw$-convergent to $x\in E$ and
$(y_{\beta})\subseteq F$ is $uaw$-convergent to $y\in F$, then
\[
x_{\alpha}\otimes y_{\beta}\xrightarrow{uaw} x\otimes y
\quad \text{in } E\widehat{\otimes}F .
\]
\end{theorem}

\begin{proof}
Assume without loss of generality that $E'$ is separable.
By Lemma~\ref{0003}, the $uaw$-topology on $E$ is metrizable by a
translation-invariant metric $d$.

Let $\widehat{E}$ and $\widehat{F}$ denote the topological completions of
$E$ and $F$ with respect to their $uaw$-topologies.
Since $d$ is translation-invariant, its completion $\widehat{d}$
metrizes the $uaw$-topology on $\widehat{E}$, and
$(\widehat{E},\widehat{d})$ is a metrizable locally solid vector lattice.

Observe that the argument of \cite[Theorem~3.2]{Den:17} remains valid for
$(\widehat{E},\widehat{d})$.
Indeed, it suffices to replace the norm $\|\cdot\|$ in the proof by the
metric $\widehat{d}$ and repeat the argument verbatim.
Consequently, there exists a subnet $(x_{\alpha_k})$ and a disjoint sequence
$(d_k)\subseteq \widehat{E}$ such that
\[
x_{\alpha_k}-x-d_k \xrightarrow{aw} 0 .
,\]
in which, $aw$ means absolutely weak convergence.
Proceeding as in the proof of Theorem~\ref{2023}, and using \cite[Corollary 6]{Z:26}
to pass to tensor products, we obtain
\[
x_{\alpha}\otimes y_{\beta}\xrightarrow{uaw} x\otimes y
\quad \text{in } \widehat{E}\widehat{\otimes}F =E\widehat{\otimes}F.
\]
This would complete 
the proof.
\end{proof}

\end{document}